\documentclass[twoside,a4paper,12pt,centertags,reqno]{amsart} 
\usepackage{amsmath,amssymb,verbatim,vmargin}
\usepackage{color}
\usepackage{tikz}
\usepackage{hyperref}
\hypersetup{colorlinks,bookmarks=true,linktocpage=true,citecolor=blue,linkcolor=magenta}

\usepackage{eulervm}   

\allowdisplaybreaks 
\usepackage{color}
\usepackage[all]{xy} 

\usepackage{mathtools}
\usepackage{MnSymbol} 

\theoremstyle{plain}

\theoremstyle{definition}

\usepackage{enumerate}

\newcommand{\bRn}{\mathbb{R}^n}

\newcommand{\cP}{{\mathcal P}}

\def\barint_#1{\mathchoice
            {\mathop{\vrule width 6pt
height 3 pt depth -2.5pt
                    \kern -9.5pt
\intop \kern -4pt}\nolimits_{#1}}%
            {\mathop{\vrule width 5pt height
3 pt depth -2.6pt
                    \kern -6.5pt
\intop \kern -4pt}\nolimits_{#1}}%
            {\mathop{\vrule width 5pt height
3 pt depth -2.6pt
                    \kern -6pt
\intop \kern -4pt}\nolimits_{#1}}%
            {\mathop{\vrule width 5pt height
3 pt depth -2.6pt
          \kern -6pt \intop \kern -4pt}\nolimits_{#1}}}
          
           \def\bariint_#1{\mathchoice
            {\mathop{\vrule width 15pt
height 3 pt depth -2.5pt
                    \kern -15.8pt
\intop \kern -8pt\intop \kern -4pt}\nolimits_{#1}}%
            {\mathop{\vrule width 9pt height
3 pt depth -2.6pt
                    \kern -10.5pt
\intop \kern -8pt\intop \kern -4pt}\nolimits_{#1}}%
            {\mathop{\vrule width 9pt height
3 pt depth -2.6pt
                    \kern -10pt
\intop \kern -8pt\intop \kern -4pt}\nolimits_{#1}}%
            {\mathop{\vrule width 9pt height
3 pt depth -2.6pt
          \kern -8pt \intop \kern -10pt\intop \kern -4pt}
      \nolimits_{  #1}}}

\def\barintlim_#1{\mathchoice
            {\mathop{\vrule width 6pt
height 3 pt depth -2.5pt
                    \kern -8.8pt
\intop \kern -4pt}\limits_{#1}}%
            {\mathop{\vrule width 5pt height
3 pt depth -2.6pt
                    \kern -6.5pt
\intop \kern -4pt}\limits_{#1}}%
            {\mathop{\vrule width 5pt height
3 pt depth -2.6pt
                    \kern -6pt
\intop \kern -4pt}\limits_{#1}}%
            {\mathop{\vrule width 5pt height
3 pt depth -2.6pt
          \kern -6pt \intop \kern -4pt}\limits_{#1}}}
          
           \def\bariintlim_#1{\mathchoice
            {\mathop{\vrule width 15pt
height 3 pt depth -2.5pt
                    \kern -15.8pt
\intop \kern -8pt\intop \kern -4pt}\limits_{#1}}%
            {\mathop{\vrule width 9pt height
3 pt depth -2.6pt
                    \kern -10.5pt
\intop \kern -8pt\intop \kern -4pt}\limits_{#1}}%
            {\mathop{\vrule width 9pt height
3 pt depth -2.6pt
                    \kern -10pt
\intop \kern -8pt\intop \kern -4pt}\limits_{#1}}%
            {\mathop{\vrule width 9pt height
3 pt depth -2.6pt
          \kern -8pt \intop \kern -10pt\intop \kern -4pt}
      \limits_{  #1}}}
          
\renewcommand{\iint}{\int \kern -3pt\int}       

\numberwithin{equation}{section}
\usepackage[symbol]{footmisc}

\makeatletter
\@namedef{subjclassname@2020}{\textup{2020} Mathematics Subject Classification}
\makeatother

\title{A one-line proof of a minimax theorem (for a notion of curvature on finite graphs)}

\author{Yi C. Huang} 
\address{Yunnan Key Laboratory of Modern Analytical Mathematics and Applications, Yunnan Normal University, Kunming 650500, People's Republic of China}
\address{Department of Mathematical Sciences, Tsinghua University, Beijing 100084, People's Republic of China}
\address{School of Mathematical Sciences, Nanjing Normal University, Nanjing 210023, People's Republic of China}
\email{Yi.Huang.Analysis@gmail.com}
\urladdr{https://orcid.org/0000-0002-1297-7674}

\date{\today} 

\thanks{Research of the author is partially supported by the National NSF grant of China (no. 11801274), 
the Visiting Scholar Program from the Department of Mathematical Sciences of Tsinghua University,
and the Open Project from Yunnan Normal University (no. YNNUMA2403).
YCH thanks Professor Peter Wakker for his indirect inspiration and helpful encouragements.}

\begin{document}

\begin{abstract}
A one-line proof of a minimax theorem due to Steinerberger is given.
\end{abstract}

\maketitle



Let $V=\{v_1, v_2,\cdots, v_n\}$ be the vertex set, $D=(d(v_i,v_j))_{i,j=1}^n$ be the distance matrix,
and $\mathbf{1}\in \bRn$ be the vector containing all 1's.
Suppose that $w\in\bRn$ satisfies
$$Dw=n\cdot\mathbf{1}.$$  
Steinerberger \cite{Ste23} views $w$ as a notion of curvature. 
Under $\min_i w_i\geq0$, he proved also a minimax theorem (``strongest" in \cite{Ste23}): $\forall$ probability measure $\cP$ on $V$
$$A:=\min_{u\in V}\sum_{v\in V}d(u,v)\cP(v)\leq\frac{n}{\|w\|_{\ell^1}}\leq \max_{u\in V}\sum_{v\in V}d(u,v)\cP(v)=:B.$$
His proof uses the von Neumann minimax theorem.
In this note we point out Steinerberger's assertion can be observed via a chain of (in)equalities in one line:
$$n=(n\cdot\mathbf{1}, \cP)=(Dw, \cP)=(w,D\cP)\in [A\|w\|_{\ell^1}, B\|w\|_{\ell^1}].$$
Proof of upper bound also goes beyond ``non-negatively curved" (i.e., $\min_i w_i\geq0$).


\section*{\textbf{Compliance with ethical standards}}


\textbf{Conflict of interest} The author has no known competing financial interests
or personal relationships that could have appeared to influence this reported work.


\textbf{Availability of data and material} Not applicable.


\bibliographystyle{alpha}

\bibliography{Hua-SteinerMinMax} 
 
\end{document}